\documentclass[11pt] {article}      %

\usepackage{amsmath,amssymb, amsfonts, amsthm, eucal}  
\usepackage{latexsym, graphicx, amscd}

\usepackage{newtxtext,newtxmath} 
\usepackage[T1]{fontenc}     
\usepackage{enumerate}       
\usepackage{float}                 

\theoremstyle{plain}
\newtheorem{theorem}{Theorem}		
\newtheorem*{corollary}{Corollary}
\newtheorem{lemma}[theorem]{Lemma}
\newtheorem{proposition}[theorem]{Proposition}   

\theoremstyle{definition}

\newtheorem{definition}[theorem]{Definition}    
\newtheorem{remark}[theorem]{Remark}		

\newtheorem{problem}{Problem}

\numberwithin{equation}{section}
\numberwithin{theorem}{section}

\def\Chron{\hbox{\rm Chron}\,}        \def\Prod{\hbox{\rm Prod}\,}
\def\Morph{\hbox{\rm Morph}\,}        \def\Ext{\hbox{\rm Ext}\,}
\def\Rais{\hbox{\rm Rais}\,}          \def\Part{\hbox{\rm Part}\,}
\def\Span{\hbox{\rm Span}\,}

\def\om{{\scriptstyle\circ}}
\def\om{\circ}
\def\e{\varepsilon}
\def\f{\varphi}
\font\caps=cmcsc10
\def\endproof{\hfill${ \vcenter{{\hrule height.4pt
    \hbox{\vrule width .4pt height 7pt\kern 7pt\vrule width.4pt}
    \hrule height.4pt}}}$}
\usepackage{titlesec}                                                           
\titlelabel{\thetitle.\ }                                                         
\titleformat*{\section}{\fontsize{14pt}{14pt} \bf}        
\usepackage{fancyhdr}
\begin{document}

\title{Induced orders in free monoids of words}

\author{Jerzy Kocik\\
{\small Department of Mathematics, Southern Illinois University, Carbondale, IL 62901}\\
{small jkocik@siu.edu}}

\date{}

\sloppy

\maketitle

\begin{abstract} 
A family of partial orders in the free monoid $A^*$ of words, induced
from a partial order in alphabet $A$, is presented.
The induced orders generalize the chronological posets that have been
defined for the two-letter alphabet only, and the morphological order.
We show that the induced orders are natural with respect to alphabet
homomorphisms.
\\[3pt]
{\bf Keywords:} Order, alphabet, words, tense systems, 
\\[3pt]
{\bf MSC:} 06A06,  
               68R15.  
\\[3pt]
{\bf Note:} \ By an {\it order} we mean a {\it partial order}.
Relations of partial order such as $\prec$, $<$, or $\subset$ are always
assumed reflexive.
\end{abstract}

\section{Motivation}  \label{s:1}

{\bf Algebra of orientons} $A^*$ is defined in [Ko] as the free monoid of
words [Lo] over a two-letter alphabet $A= \{\,\pi , \varphi\,\}$ ordered $\pi\prec\varphi$.
Two partial orders [Be] are introduced into the monoid $A^*$:
the self-evident {\bf morphological} order representing complexity
of words, and a less apparent {\bf chronological} order that has been
induced from $\pi\prec\varphi$ assumed in the alphabet $A$.
The relations are defined as follows:

\begin{definition} 
Let $A^*=\{\pi ,\varphi\}^*$ be words
 over $A=\{\pi,\varphi\}$.  For any two elements $v,w\in A^*$ one defines

\begin{enumerate}
\item[(a)]
morphological order:  
$v<w$  if $w$ may be obtained
by inserting some letters of $A$ into $v$.
\item[(b)] 
chronological order:  $v\prec w$ if $w$ may be obtained from $v$
by erasing some letters $\pi$ and/or inserting some letters $\varphi$ into $v$.
\end{enumerate}
%
%
\end{definition}

\def\ignorea{
\begin{definition} 
Let $A^*=\{\pi ,\varphi\}^*$ be words
 over $A=\{\pi,\varphi\}$.  For any two elements $v,w\in A^*$ one defines
\\
  a) morphological order:
\\
  \hskip3em ~$v<w$ \ \  if \ $w$ may be obtained
                       by inserting some letters of $A$ into $v$.
  b) chronological order:
$$
 v\prec w \quad \hbox{if}\qquad 
   \left\{ 
  \begin{array}{cl}
   &w \hbox{ may be obtained from}\ v\ \hbox{by erasing some letters}\ \pi \cr
   &  \hbox{ and/or inserting some letters}\ \varphi \hbox{into}\ v.
       \end{array}\right.
$$
\end{definition} 
}

If $<\!\!<$ and $\prec\!\!\prec$ denote the covering relation
(immediate succession) of the respective orders, then the above
definitions may be expressed in terms of single insertions:
\begin{enumerate}
\item[(a')]
$v<\!\!< w$ \ \ if \ \ $\exists x\in A :\
  v'\om x\om v''=w$ \ \ for some splitting $v=v'\om v''$
\item[(b')]
$v\prec\!\!\prec w$ \ \ if \ \
  either $v'\om\varphi\om v''=w$ \ or \ $v=w'\om\pi\om w''$
for some splitting $v=v'\om v''$ \ or \ $w=w'\om w''$.
\end{enumerate}

The corresponding posets are denoted respectively $\Morph A^* =
\{\,A^*, <\,\}$  and $\Chron A^* = \{\,A^* , \prec\,\}$.
Figure 1 displays them for words of length $|w|<4$.
Poset $\Morph A^*$ has the least element, the empty word $\e$;
$\Chron A^*$ is unbounded.

\begin{figure}[h]  
\[
\includegraphics[width=4in]{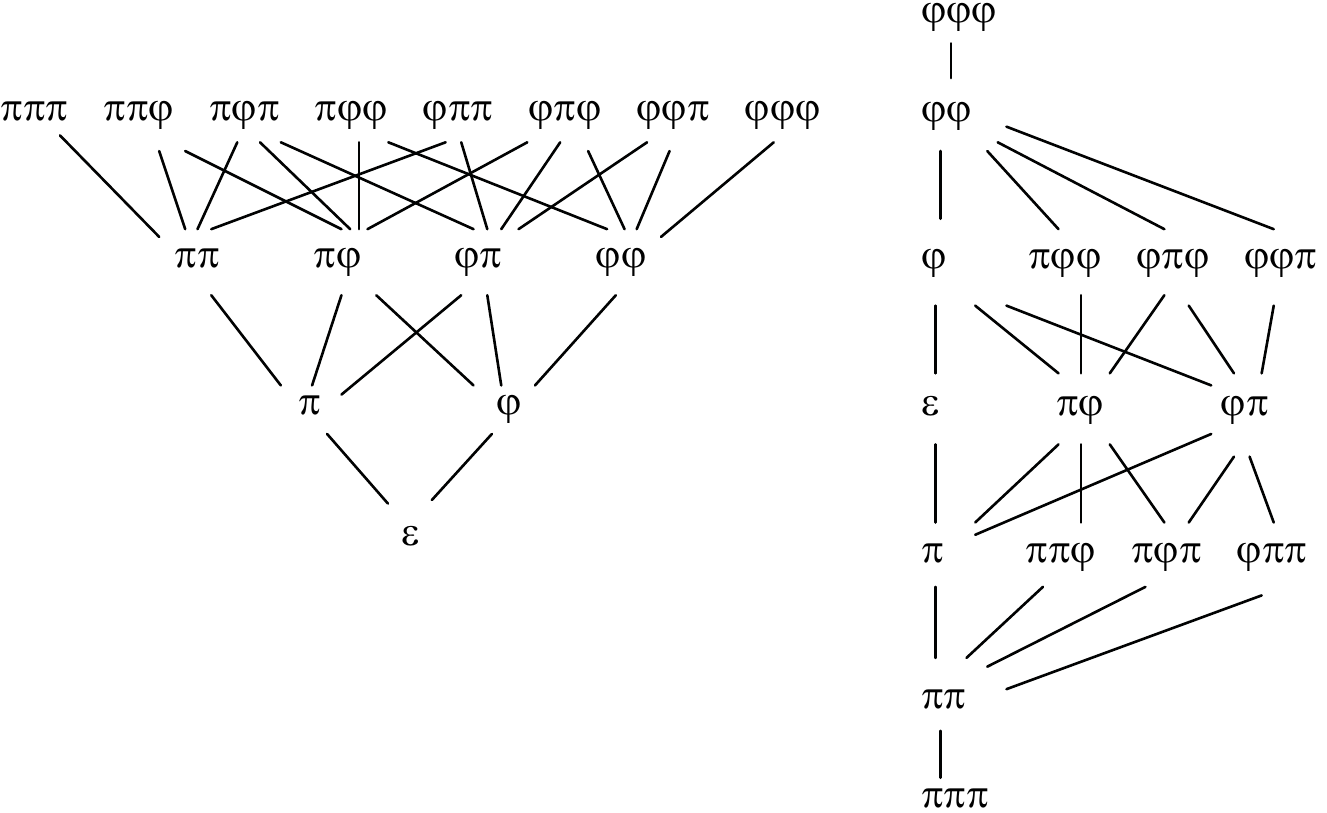}
\]
\caption{The morphological and chronological orders of $A^*=\{\pi,\varphi\}$.}
\label{fig:1}
\end{figure}

\begin{remark}
Algebra of ``orientons" $\{\pi,\varphi\}^*$ was introduced in [Ko] to model
chronological meanings of grammatical tenses.
Intuitively, $\pi$ refers to the ``past", $\varphi$ --- to the ``future",
the empty word $\e$ --- to  the ``present tense," and,
for instance, $\pi\pi\varphi$ --- to the ``future-in-the-past-in-the-past".
The order $\pi\prec\varphi$ in the alphabet (``past" precedes ``future")
induces a ``chronological order" among the words of $\{\pi,\varphi\}^*$.  This
infinite set contains a homomorphic images of formal tense systems
including a mathematical model of the tense system of English.
\end{remark}

While the morphological order is a natural relation associated with any
free monoid over an alphabet, the chronological order in $\{\pi,\varphi\}^*$
is not, and the question rises, ``How far may this construction be
extended beyond the simple two-letter alphabet?"  The form of Definition~1
hardly seems to suggest any possible natural generalization.
\\

\begin{problem}   
Let $A$ be a partially ordered set of a countable cardinality,
and let $A^*$ be the free monoid of words over the alphabet $A$.
Is there any natural order in $A^*$ which
\  (i) would be an extension of the order in $A$ \
($A$ is embedded into $A^*$ as the one-letter words);
\  (ii) would be natural with respect to homomorphisms of alphabets; and
\  (iii) would coincide with the chronological order for the simple case of
$A=\{\pi\prec\varphi\}$ ?
\end{problem}

The answer to this problem is affirmative, although the construction
is not a direct generalization of Definition~1.1b.
The following notes describe this construction.

\section{Induced order}  \label{s:1}

In this section, we show a {\it principal} construction by which
every letter of an ordered alphabet defines an induced order among the
words over the original alphabet without the chosen letter.

Let $\{\,A,\prec\,\}$ be a partially ordered set.  For any $a\in A$ we
denote $A_a = A - \{a\}$.  Consequently,
$$
   A^*_a = \big( \, A - \{a\}\, \big)^*
\eqno(1)
$$
denotes the free monoid over the alphabet $A_a$.
Of course, the relation $\prec$ restricted to $A_a$ makes it a poset.
We shall call a word $w' \in A^*$ an $a$-{\it extension} of a word
$w\in A^*_a$ if $w'$ may be obtained by a number of insertions of the
letter $a$ into $w$.
(For instance {\it mississippi} is an $s$-extension of {\it miiippi},
which in turn is a $p$-extension of {\it miiii}).

Now, we define a relation among words of $A_a^*$, induced by the order
$\prec$ in $A$.

\begin{definition} 
Let $v,w \in A^*_a$ be two words over a partially ordered alphabet
$\{ A_a,\prec\}$.
We say write
$$
                    v \prec_a w
$$
if $a$-extensions $v'$ and $w'$ of $v$ and $w$ respectively exist, such
that they are of the same length $|v'|=|w'|=n$,  and  $v'\prec w'$ in
$A^n$, i.e., for each letter it is $v'_i \prec w'_i$ in $A$, \
$i=1,\ldots,n$.
\end{definition} %

It is not clear whether different insertions would not lead to different
directions of $\prec_a$ for the same pair of words.
Here is the main assertion of this note:

\begin{theorem} 
The relation $\prec_a$ defines a partial order in $A^*_a$.
\end{theorem} %

The relation $\prec_a$ is clearly reflexive and transitive.
The problem is whether it is skew-symmetric. First, we prove the
following lemma.

\begin{lemma} 
Let $w'$ and $w''$ be two $a$-extensions of the same length
$|w'|=|w''|=n$ of a word $w\in A^*_a$ such that $w'\prec w''$ in $A^n$.
Then $w' = w''$.
\end{lemma} 

\begin{proof}   
Let \ $w' = w'_1\om w'_2\om\ldots\om w'_n$,
where $w'_i\in A$ for each $i=1,\ldots ,n$.
Since $w''$ is of the same length as $w'$, it must be composed
from the same collection of letters, and the arrangement of letters in
$w''$ is a permutation of the arrangement of letters in $w'$:
$$
   \qquad\qquad w''=w'_{\sigma(1)}\om w'_{\sigma(2)}\om\ldots\om w'_{\sigma(n)}
                 \qquad\hbox{ for some} \ \sigma\in S_n.
$$
Now, by definition the assumed relation $w'\prec w''$ in $A^n$ means
that it holds for the corresponding letters of $w'$ and $w''$ in $A$:
$$
\begin{array}{ccccc}
      w'_1 &\prec& w''_1 &=& w'_{\sigma(1)}\cr
      w'_2 &\prec& w''_2 &=& w'_{\sigma(2)}\cr
     \vdots&     &       & & \vdots\cr
      w'_n &\prec& w''_1 &=& w'_{\sigma(n)}
\end{array}
$$
Each permutation may be uniquely decomposed into a number of cycles.
Assume that $\sigma$ has a cycle of order bigger than 1, say $k$.
Then, for some $w'_i\in A$ it is
$$
\begin{array}{rcl}
   w'_i    \quad\quad   &\prec\quad &w'_{\sigma(i)}    \cr
   w'_{\sigma(i)}  \quad    &\prec\quad &w'_{\sigma(\sigma(i))} \cr
    \ldots      \quad    &\quad      &\ldots  \cr
   w'_{\sigma^{k-1}(i)}     &\prec\quad &w'_{\sigma^k(i)} = w'_i
\end{array}
$$
By transitivity both $w'_i\prec w'_{\sigma (i)}$ and
$w'_{\sigma (i)}\prec w'_i$ hold.  This contradicts  the partial order
of the alphabet,  unless the order of each cycle of $\sigma$ is 1.
Therefore $w'_i = w'_{\sigma(i)}$ for each $i$, so $w'=w''$,
proving the lemma.
\end{proof}

\begin{corollary} 
Any two $a$-extensions of the same length n of a word $w\in A^*_a$ are
either identical or incomparable in the product poset $A^n$.
\end{corollary}

Now we can prove the theorem.
\\

\noindent
{\bf Proof of Theorem 2.2:}
In order to show that the induced relation $\prec_a$ in $A^*_a$
is skew-symmetric for different elements,
let us assume {\it a contrario} that there exist
two different insertions of the letter $a$ into a pair of words $v$ and
$w$ in $A^*_a$,  such that the resulting two pairs of $a$-extensions,
$v'$, $w'$, and $v''$, $w''$, lead to opposite relations:
\\

\hskip.8in $v\prec_a w$ \ \ by one insertion, \
                  for which $v'\prec w'$ in $A^n$, \ and\hfill
\\[-7pt]

\hskip.8in $v\succ_a w$ \ \ by the other one, \
                  for which $v''\succ w''$ in $A^k$,\hfill
\\[7pt]
where the lengths of the words are $|v'|=|w'|=n$ and $|v''|=|w''|=k$
for some $n$ and $k$.  In the form of a diagram:
%
$$
\begin{array}{ccccccc}
    &        & v  &\!\!\!\!\leftrightarrow\!\!\!\!          & w \hfill  \cr\cr
    &\diagup &    &\diagup\!\!\!\!\!\diagdown &   &\diagdown \cr\cr
 v'\!\! &\prec   & w'\!\! &                         &\!\! v'' &\succ &\!\! w''\cr\cr
  (\hbox{in} & A^n)&            &&   & (\hbox{in} & A^k)   
\end{array}
$$

Assume that $v\neq w$, since Lemma 2.3 proves the theorem for $v=w$.

Notice that additional simultaneous insertion of the letter $a$ into the
words $v',w'$ at the {\it same} position preserves the original
relation, now in $A^{n+1}$. A number of such insertions will be called
a {\it coherent} $a$-extension of a pair of words (of the same length).

Now, since $w'$ and $w''$ result by insertions of the letter $a$
into the same word $w$, one may find a further minimal
$a$-extensions of $w'$ and of $w''$ such that the resulting
words $\bar w'$ and $\bar w''$ will be identical
in $A^m$ for some $m$:
$$
           \bar w'\quad=\quad\bar w''    \qquad\hbox{in~~}A^m
\eqno(2)
$$
If $\bar v'$ is the $a$-extension of $v'$ coherent with that of $w'$,
and $\bar v''$ the extension coherent with that of $w''$,
the relations between the words become:
$$
\begin{array}{rl}
      \bar v'\qquad\prec\qquad&\bar w'\cr
                              &\bar w''\qquad\prec\qquad\bar v'',
\end{array}
\eqno(3)
$$
in $A^m$, and by (2) and transitivity:
$$
              \bar v'\qquad\prec\qquad\bar v''.
\eqno(4)
$$
By Lemma 2.3, this implies \ $\bar v' = \bar v''$.
Therefore, from (3):
$$
      \bar v'\quad = \quad\bar w'\quad  = \quad\bar w''\quad
                   = \quad\bar v'' \eqno{\hbox{(in~~} A^m) \qquad}
$$
Removing all $a$'s from these words we  get $v = w$,
which contradicts the assumption and concludes the proof.
\endproof
\\

We obtain a whole family of induced orders, labeled by the elements of $A$.
In particular:

\begin{corollary} 
If $A=\{\,\pi ,\eta ,\varphi\,\}$ is linearly ordered $\pi\prec\eta\prec\varphi$,
then the induced order $\prec_\eta$ in $A^*_\eta=\{\pi ,\varphi\}^*$ coincides
with the chronological order of the word algebra over the alphabet
$\{\pi ,\varphi\}$.    (see Definition~1.1).
\end{corollary} %

\section{Augmentation}  \label{s:3}

Now the solution to the Problem (Section 1) seems plausible.
In order to obtain a relation in the monoid $A^*$ over an ordered
alphabet $A$,
one has to enrich first the alphabet by one element, say $e$, and to
extend the order of $A$ into $A'=A\cup\{ e\}$, and then apply the
technique of induced order described in the previous section.
We shall call letter $e$ an {\it auxiliary\/} letter.
Poset $A'$ will be called an {\it augmented alphabet}.

\begin{definition} 
An augmentation of a (finite) poset $A$ is an isomorphism of $A$
into a poset $A'$ of cardinality $|A'|=|A|+1$.
\end{definition}

The original alphabet is restored by dropping letter $e$, i.e.
as a set $A\equiv(A\cup\{e\})_e$.
The partial order $\prec_e$  defined by Definition~1.1  turns $A^*
\equiv (A\cup\{ e\})^*_e$, due to Theorem 2.3, into a poset.

Of course, the order so obtained strongly depends on the particular
choice of augmentation.
\\

\noindent
{\caps Example} 1: \
Consider $A=\{\pi ,\varphi\}$ with $\pi\prec\varphi$. In order to get an induced
order in $A^*$, an augmented poset must be constructed with elements
$A'=\{\pi ,\varphi ,\eta\}$, where $\eta$ is an auxiliary letter.
There are three ways to equip $A'$ with a {\it linear} order that
agrees with the order in $A$:
$$
     (a) \qquad \eta\prec\pi\prec\varphi, \qquad\qquad
     (b) \qquad \pi\prec\eta\prec\varphi, \qquad\qquad
     (c) \qquad \pi\prec\varphi\prec\eta
$$
Each leads to another partial order in $A^*=\{\pi,\varphi\}^*$. Figure~2
displays the corresponding induced posets for the words of length
$|w|\leq 3$.
Case (b) is identical with the chronological order (see Definition~1.1).
Case (c) is dual to the case (a) \ (replace $\varphi$ with $\pi$
and flip the diagram upside down).
Notice that in each of these cases, the one-letter words, which
may be identified with the elements of the alphabet,
preserve their order $\pi\prec\varphi$ within $A^*$.

\begin{figure}[H]  
\[
\includegraphics[width=4.8in]{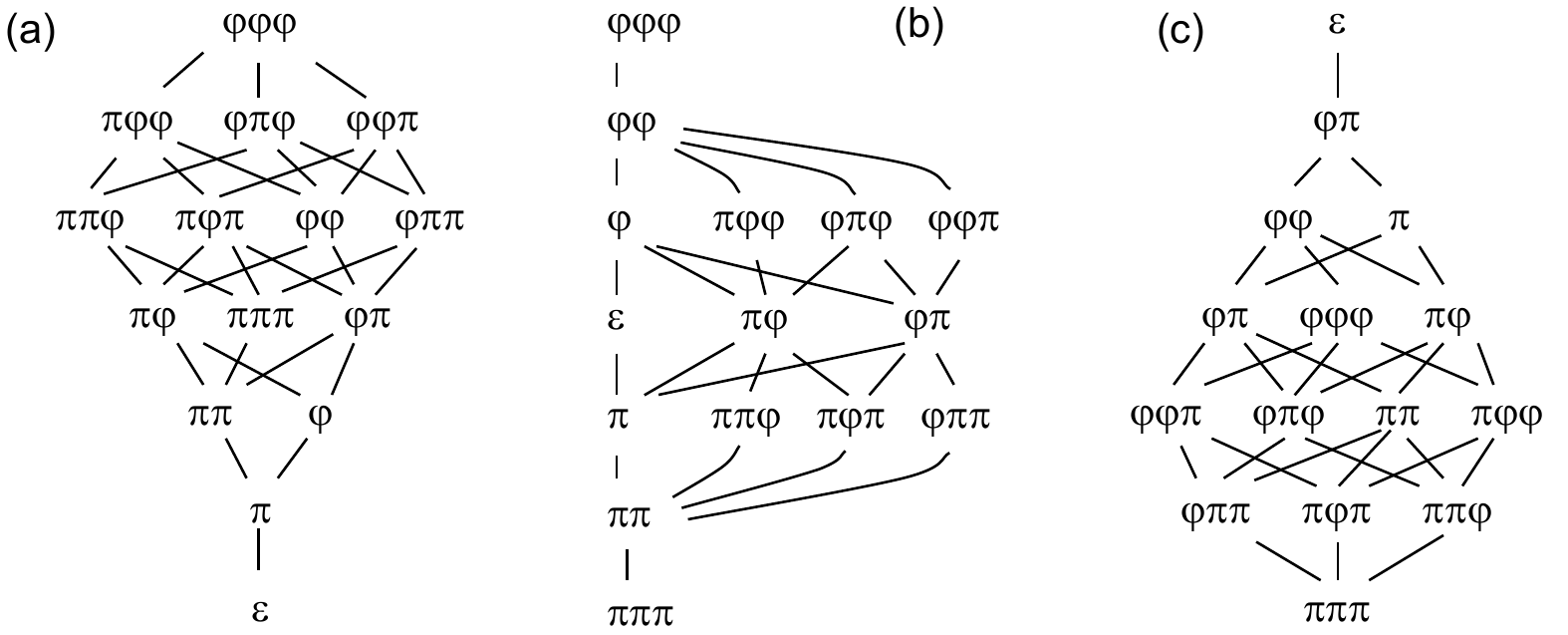}
\]
\caption{Three posets over $\{\pi ,\varphi\}^*$ induced by 
         linear augmentations of $\pi\prec\varphi$.}
\label{fig:2}
\end{figure}

The last observation can be generalized:

\begin{proposition} 
  For any $n\in\mathbb N$ and for any augmentation of a poset $A$,
  the poset $A^n$ with the product partial order is isomorphically
  embedded into the induced word posets $A^*$.
\end{proposition} %

\begin{proof}   
Identify the Cartesian product $A^n$ with words in $A^*$
of the fixed length $n$, $A^n \equiv \{ a\in A^*\bigm| |a|=n\,\}$.
For any augmentation, if two words are related in $A^n$,
so are they, by definition, in $A^*$; if they are not related in
$A^n$, then by Corollary~2.4 they are not related in $A^*$.
\end{proof}

In particular, the natural embedding of an alphabet $A$ into the
one-letter words in $A^*$ is an isomorphism of the order structures.
For an illustration of  $n=2$ and $n=3$, recognize the particular posets
of $A^n$ (Figure~4) in Figure~2 and Figure~3.
Notice that the above property may be extended to an embedding of
$A'=A\cup\{e\}\to A^*$, if the empty word $\e\in A^*$  is reinterpreted
as the auxiliary letter $\eta$ in $A'$.

The following obvious property ensures naturality of induced order,
which was sought in Problem~1.2.

\begin{proposition} 
   Let $f:\  A\to B$ be a homomorphism of posets.
   For any $e\in A$, the induced map  $f^*:\  A^*_e\to B^*_{f(e)}$,
   defined letter-wise, is also a poset homomorphism of induced orders.
\end{proposition} %

\begin{proof}   
Proof is straightforward. Let $f\colon A\to B$ be a poset
homomorphism, i.e. if $a\prec b$ in $A$, then $f(a) \prec f(b)$ in $B$.
Let $f^*\colon A^*\to B^*$ be a letter-wise extension of $f$.
Clearly, it may be restricted to $f^*\colon A^*_e\to B^*_{f(e)}$.
Relation $v\prec_e w$ in $A^*_e$ means that there are $e$-extensions
of $v$ and $w$, such that $v'\prec_e w'$ in $A^n_e$ for some 
$n\in\mathbb N$, i.e. $v'_i\prec_e w'_i$ in $A_e$ \ for \ $i=1,\ldots,n$.
So, $f(v'_i)\prec_{f(e)} f(w'_i)$ in $B_e$, and therefore
$f^*(v')\prec_{f(e)} f^*(w')$ in $B^n_e$.
Hence, by definition, $f^*(v)\prec_{f(e)} f^*(w)$ in $B^*$.
\end{proof}

In particular,

\begin{corollary} 
(i) \ If $B\subset A$ as sets, and the partial order of $B$ is
that of $A$ restricted to $B$, then the induced poset $\{\,B^*_e,\prec_e\,\}$
is a subposet of $\{\, A^*_e , \prec_e\,\}$, for any $e\in B$.
\\
(ii) If $\prec'$ be a suborder of a partial order $\prec$ in $A$, then for
any $a\in A$ the order $\prec'_a$ is a suborder of $\prec_a$ in $A^*_a$.
\end{corollary} %

As an extremely simple illustration of (i), consider $\{\varphi\}$ as a
one-element subposet of $A'=\{\pi,\eta,\varphi\}$ for any of the given examples
of augmentation. The word algebra $\{\varphi\}^*$ consists of powers $\varphi^n$.
In a particular example either $\varphi\succ\eta$ or $\varphi\not\sim\eta$,
and hence $\{\varphi\}^*$ is either linearly ordered or trivial, and so it
occurs in the corresponding posets $\{\pi,\varphi\}^*$.
For an illustration of (ii), compare Figure~3 with Figure~2, where the
corresponding Hasse diagrams form subgraphs on the alphabet level, as
well as in the word algebras.
\\

\noindent
{\caps Example} 2: \
Consider two {\it nonlinear} extensions of $\pi\prec\varphi$:
$$
\begin{array}{cccccc}
      &\varphi&         &      &       & \eta  \\
      &       &\!\!\diagdown\!\!&      &\!\!\diagup\!\!\\ 
  (a) &       &         &\pi   & 
\end{array}
\qquad
\begin{array}{cccccc}
      &         & \varphi&          &    \cr
      &\!\!\diagup\!\! &        &\!\!\diagdown\!\!      \\ 
  \pi &         &        &          &\eta      & (b)     
\end{array}
$$
These lead to posets, which are displayed in Figure~3
for words $|w|\leq3$.

\begin{figure}[H]  
\[
\includegraphics[width=4in]{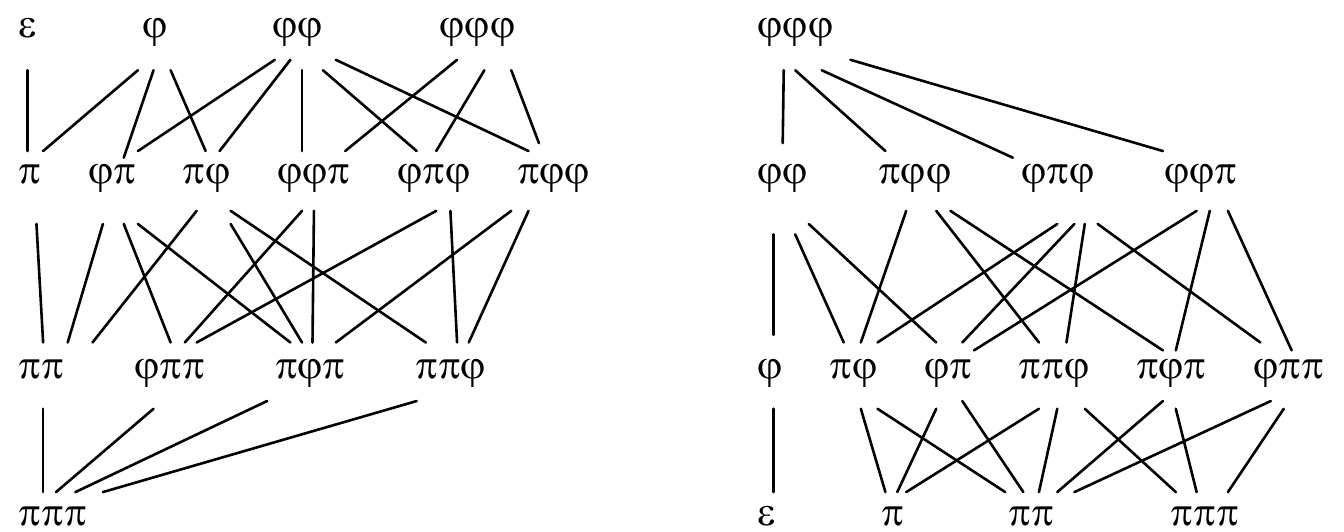}
\]
\caption{Orders in $\{\pi,\varphi\}^*$ induced from non-linear augmentations.}
\label{fig:3}
\end{figure}

Although posets of Example 2 look ``strange," notice that each is a
particular suborder  of two of the posets considered in Example~1.
This is because the orders of $A$ in Example~1 are particular
{\it linearizations}  of the alphabet posets (a) and (b) above).

\section{Further examples and applications}  \label{s:4}

Now let us review a few special cases, illustrated by rather simple examples.

\begin{definition} 
A {\it raising} augmentation of a poset A is a poset $A' = A\cup\{ e\}$
with the partial order this of $A$ complemented by relation $e\prec a$\
for any $a\in A$.  The induced poset will be denoted $\Rais A^*$.
\end{definition} %

For illustration of $\Rais\{\pi\prec\varphi\}^*$ see Example~1a.

Augmentation may also be applied to mere sets (viewed as posets with
the trivial order). In particular:

\begin{corollary} 
The induced order of a trivial poset is the morphological order.
\end{corollary} 

\noindent
{\caps Example} 3: \
Let $A=\{\varphi ,\pi\}$ be a set. Raising augmentation of $A$
into a poset $A'=\{\pi,\varphi,\eta\}$ with a two-step relation:
$$
\begin{array}{ccccc}
    \varphi&                 &     &                & \pi  \\
           &\!\!\diagdown\!\!&     &\!\!\diagup\!\!  \\ 
           &                 &\eta & 
\end{array}
$$
results in morphological order of orientons. (See Definition~1.1 and Figure~1.)
Quite surprisingly, both key orders of the `algebra of orientons'
are describable in terms of induced orders.

\begin{definition} 
A trivial augmentation of a poset $A$ is a poset $A'=A\cup\{e\}$
with the auxiliary letter $e$ left unrelated to $A$.
\end{definition} %

\begin{corollary} 
The poset induced from the trivial augmentation of a poset $A$ is
the disjoint sum of product orders in subsets of $A^*$:
$$
     \{\,A^*,\prec\,\}= \Prod A^0+\Prod A^1 +\ldots+\Prod A^k+\ldots
$$
where $\Prod A^k$ is the Boolean lattice of the product order among the
words of a fixed length.
(Clearly, $\Prod A^0\cong \{\e\}$, and $\Prod A^1\cong \{A,\prec\}$).
\end{corollary} %

\begin{proof}   
By Proposition~3.2, two words of the same length, $|w|=|v|=k$,
are related in $A^*$ in the same way as in $\Prod A^k$.
Words of different lengths in $A^*$ are not related:
any $e$-extensions of $w$ and $v$ resulting in the same length must have
different numbers of the letter $e$, with some of them occurring where a
letter of $A$ appears in the other word. Since $e\not\sim A$, the
extended words are incomparable, and therefore so are the original words
$w$ and $v$.
\end{proof}
\\

\noindent
{\caps Example} 4: \
Trivial augmentations of the two-letter poset $A=\{\pi\prec\varphi\}$
with auxiliary letter $\eta$ is left unrelated to $A$:
$$
\begin{array}{ccccc}
    \varphi&                 &     &                & \eta  \\
           &\!\!\diagdown\!\!&     &\!\!\diagup\!\!  \\ 
           &                 &\pi & 
\end{array}
$$
splits the word algebra $A^*$ into a family of disconnected Boolean lattices
of constant word length, as illustrated in Figure~4.

\begin{figure}[h]  
\[
\includegraphics[width=2.5in]{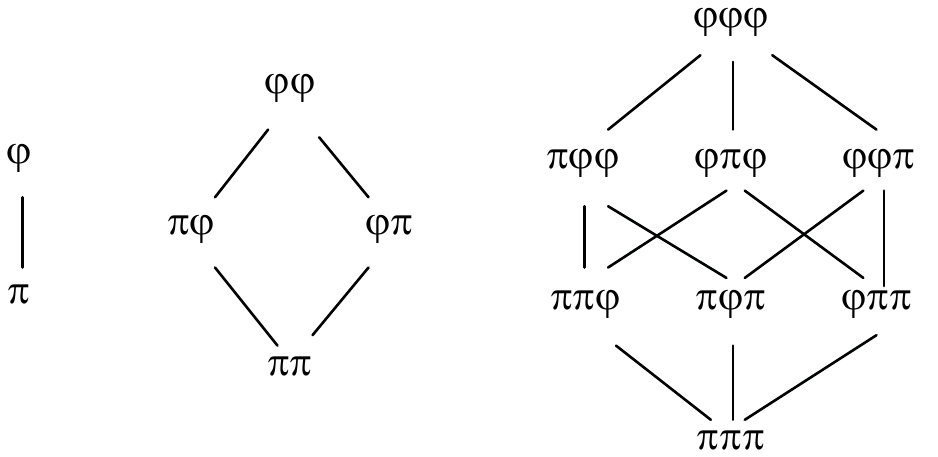}
\]
\caption{Order in $\{\pi ,\varphi\}^*$ induced from trivial augmentation.}
\label{fig:1}
\end{figure}

Note, the range of the induced orders over the same word monoid:
the words of {\it connected} pieces of $\Prod A^*$ appear as the
horizontal layers in $\Morph A^*$ \ (Cf. Figure~1 and Figure~4).
\\


\noindent
{\caps Example} 5: \
Consider the poset $A^*_\eta$ induced from the following augmentation of a
trivial poset
$$
\begin{array}{ccccc}
    \varphi&                 &     &                & \pi  \\
           &\!\!\diagdown\!\!&     &                   \\ 
           &                 &\eta & 
\end{array}
$$
(letter $\eta$ related to $\varphi$ only).
The resulting partial order is illustrated in Figure~5.

\begin{figure}[h]  
\[
\includegraphics[width=4in]{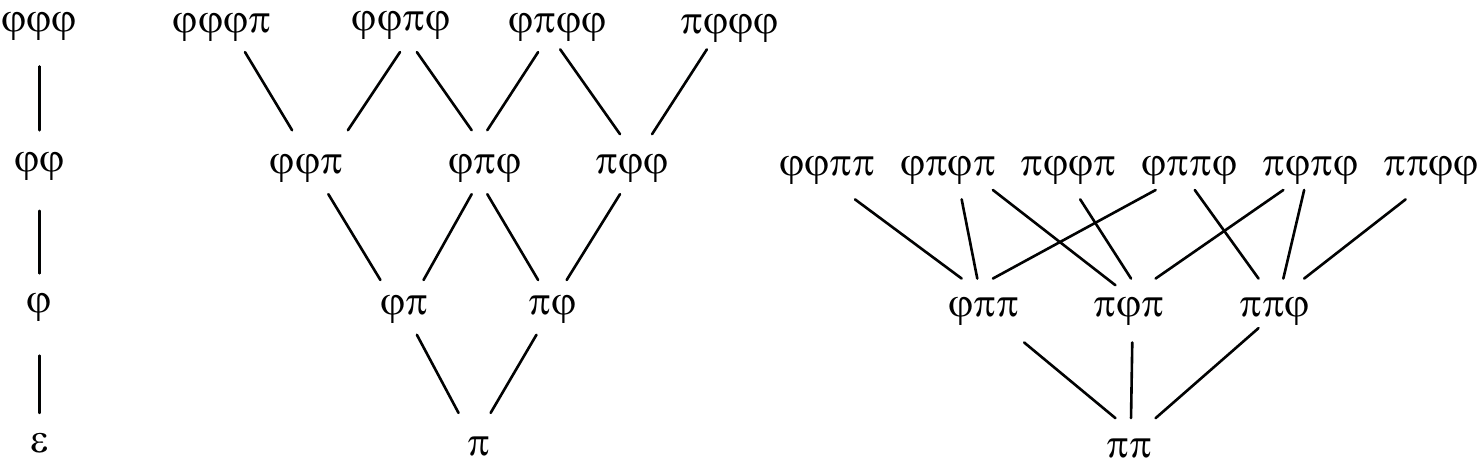}
\]
\caption{The partition poset for $A=\{\pi,\varphi\}$.}
\label{fig:5}
\end{figure}

If the letter $\pi$ is viewed as a {\it separating bar},
$``|"$, then the $i^{th}$ connected piece of the above graph shows the
possible distributions of a number of items $\varphi$ into $i$ boxes,
including the partial order of such distributions.
The $i^{\rm th}$ piece has has $\pi^{i-1}$ as the least element
($i$ empty boxes), and is isomorphic to the product $\mathbb N^i$.
This suggests the following:

\begin{definition} 
   Let $A$ be a poset. Consider an augmented poset
   $A''= A\cup\{ e ,\,\mid\,\}$  with the extended partial order:
   $e\prec\,A$, and $\mid$ unrelated to $A$ or $e$.
   A partition poset $\Part A^*$ is the word algebra over the alphabet
   $(A\cup \{\mid\})\equiv A'$ with the order induced from $A''$.
\end{definition} %

\section{Summary}  \label{s:6}

Each poset $A$ treated as an alphabet leads to a natural family of
well-defined {\it induced} partial orders in the set of
 words over this alphabet (Theorem 2.2) .
The {\it principal} construction of the induced order goes via dropping
a letter, say $a$, from the alphabet $A$, and considering the new set,
$A_a$, as the alphabet for words, $A^*_a$.  The choice of the letter to
be dropped determines the partial order in $A^*_a$.
{\it Augmentation} allows an induced order to be defined between the
words over the initial alphabet $A$ by, first, embedding the alphabet $A$
into a larger poset $A'=A\cup\{ e\}$, and then applying the
principal construction by dropping the auxiliary letter $e$.
\\

Natural properties of induced order easily follow
(expressed here for {\it augmentation}).
The induced order is an extension of the order in the alphabet~$A$
$$
         A \ \xrightarrow{\ \ \hbox{homo}\ \ } \ A^*
$$
The construction is {\it natural} with respect to the homomorphisms of
alphabets (Proposition~3.3), making the following diagram commute:
$$
\begin{array}{ccc} 
 A^* & \xrightarrow{\ \ \ f^*\ \ \ }  & B^*    \\[7pt]
   {\Big\uparrow}         &              &{\Big\uparrow}    \\[7pt]
            A   & \xrightarrow{\ \ \ f \ \ \ \ }  & B      
\end{array}
$$
where $f^*$ is a letter-wise homomorphism induced from $f$.
The property that contrasts the induced orders with
the lexicographical order is that  the product posets $\Prod A^n$
are isomorphically embedded into $A^*$ \ (Proposition~3.2).

Since the induced order {\it extends} that of the alphabetic order, let us
denote the induced poset as $\Ext A^*$ for augmentation,
or $\Ext\,A^*_e$ for the principal construction:
$$
             \Ext\,A^*_a = \{\, A^*_a,\, \prec_a\,\}
$$
A few canonical constructions (by augmentation, $e\not\in A$) may be indicated:
$$
\begin{array}{rcll}
  \Morph A^*       &=& \Ext\{e\prec A\}^*_e        &( A \hbox{\ is a set})\qquad\cr
  \Prod A^*        &=& \Ext\{e\not\sim A\}^*_e     &\cr
  \hbox{Rais}\,A^* &=& \Ext\{e\prec A\}^*_e        &\cr
  \hbox{Part}\,A^* &=& \Ext\{e\prec\,A\;,\ \eta\not\sim A\,\}^*_e
                                             &(\eta\not\in A,\ A \ \hbox{is a poset})\cr
  \Span A^*        &=& \Ext\,\{L(A)\prec e\prec G(A)\}^*_e  
\end{array}
$$
where in the last poset $L(A)$ and $G(A)$ are the least and the greatest
elements of $A$ respectively.
\\

The above augmentations may be illustrated symbolically:
$$
\includegraphics[width=4.8in]{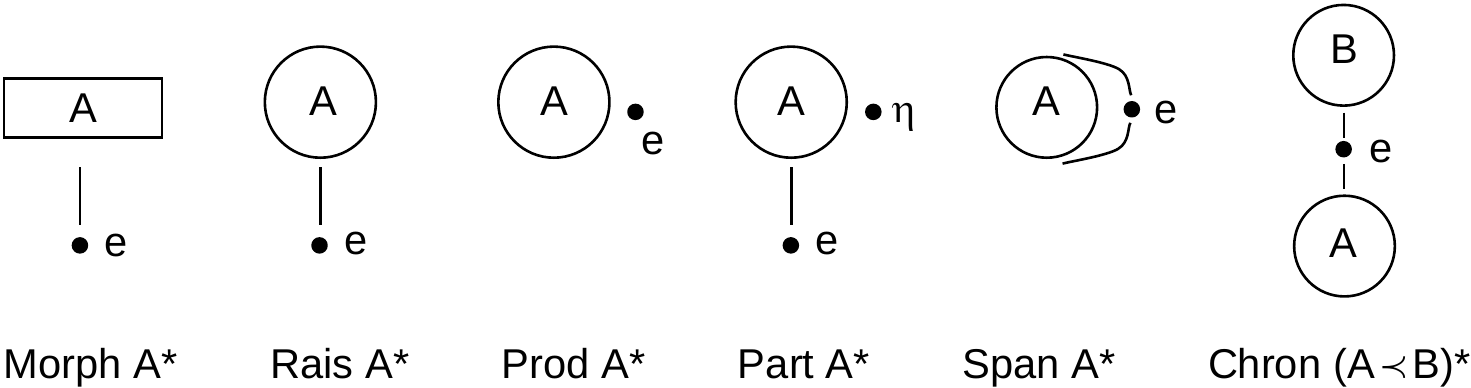}
$$

%
%

As to the algebra of orientons $A^*=\{\pi,\varphi\}^*$,
surprisingly both the chronological and the morphological order turns
out to be induced extensions:

$$
\begin{array}{rl}
\Chron A^* &= \Ext \{\pi\prec\eta\prec\varphi\}^*_\eta  \cr
\Morph A^* &= \Ext \{\pi\succ\eta\prec\varphi\}^*_\eta
\end{array}
$$

Another construction (interesting in the context of discrete models of
causal properties of space-time) concerns the union $A\cup B\cup \{e\}$
of posets $A$ and $B$, complemented by $A\prec e\prec B$.
The induced order in $(A\cup B)^*$ defines a poset
$$
        \Chron(A\prec B)^* = \Ext\{A\prec e\prec B\}^*_e
$$
which may be viewed as a direct generalization of the chronological
order.   
By analogy to relativity theory, it seems natural to define
the {\it future cone} and the {\it past cone} as the image of $A^*$ and
$B^*$ in $(A\cup B)^*$, respectively, and the {\it elsewhere} as
$(A\cup B)^* - (A^* \cup B^*)$.
\\

Here are some questions concerning the induced orders:
How are particular properties of the ordered alphabet
reflected in the induced order of words;
How does the induced order relate to the  ``algebra of products
of partial orders" (in the sense of [2]);
What is the relationship between algebra of partial orders
(treated on the level of alphabets) and that lifted to the words;
How is the structure of the alphabet $A$ reflected in the structure
of the family of posets obtained by deleting different letters from $A$.
In particular, notice ``non-commutativity" of the construction;
although $(A_a)_b^* = (A_b)_a^* = (A-\{a,b\})^*$, but
$\Ext\,(A_a)^*_b \neq\Ext\,(A_b)^*_a$.

%


\end{document}